\documentclass[reqno]{amsart}
\usepackage{amssymb}
\usepackage{t1enc}
\usepackage[latin1]{inputenc}
\usepackage[english]{babel}
\begin{document}
\title{ Elliptic functions, infinite products and modular relations}
\author[]{Abd Raouf Chouikha} \footnote
{chouikha@math.univ-paris13.fr. 4, Cour des Quesblais 35430 Saint-Pere, France   }


\begin{abstract}
Infinite products expansions of the Weierstrass elliptic function \ $\wp(z) = \wp(z,1,\tau)$\ and $n$-order transformations allow us to provide some modular relations. 
\end{abstract}
\maketitle 
\section{Introduction}

Let \ $z, \ \tau$\ be complex numbers such that \ $\Im \tau > 0$\ and \ $\Im z < 2 \ \Im \tau,$\ and 
consider the Weierstrass elliptic function \ $\wp(z) = \wp(z,1,\tau)$\ with two primitive periods $2,2\tau$.\\ 
The $n$-order odd decomposition as infinite products (see for example Tannery and Molk [10, T.2,p.246]) of the Weierstrass sigma functions allows us to deduce the $n$-decomposition of the elliptic Weierstrass function \ $\wp(z,\tau)$   
$$\wp(nz,n\tau) - e_1(n\tau) = \left(\frac{4}{\pi^4}\right)^{n-1} \frac{\theta_3^2(0,n\tau) \theta_4^2(0,n\tau)}{\left[\theta_3^2(0,\tau) \theta_4^2(0,\tau)\right]^n} \ \prod^{n-1}_{m=0}  \left[\wp(z+\frac{m}{n},\tau) - e_1(\tau)\right],$$
where \ $e_1 = \wp(1,\tau),$ as well as for its derivative
$$\wp'(nz,n\tau) = \left(\frac{4}{\pi^4}\right)^{n-1}\frac{{\theta'}_1^2(0,n\tau)}{\left[{\theta'}_1^2(0,\tau)\right]^n} \ \prod^{n-1}_{m=0}  \wp'(z+\frac{m}{n},\tau).$$

For any odd integer \ $n$\ we prove the following modular identities :
$$ \sum_{k\neq 0} \left[ \frac{1}{\sin(2kn\pi \tau+n\pi z)}\right] = \sum_{m=0}^{n-1} \sum_{k\neq 0} \left[ \frac{1}{\sin(2k\pi \tau+\pi (z+\frac{2m}{n}))}\right] =$$
$$2^{1-n}\sum_{k\neq 0} \prod_{m=1}^{n-1} \left[ \frac{1}{\sin(2k\pi \tau+\pi (z+\frac{2m}{n}))}\right] =  \frac{\theta_{2}^2(0,n\tau)}{\left[\theta_{2}^2(0,\tau)\right]^n}\ \prod^{n-1}_{m=0} \sum_{k\neq 0} \left[ \frac{1}{\sin(2k\pi \tau+\pi z+ \frac{m\pi}{n})}\right].$$

To prove that we will use the $n$-decomposition of \ $\wp(z,\tau)$\ 
as well as the $n$-transformations of:
$$ \frac{\wp'(z)}{\wp(z) - e_1} =  -\frac{2\pi}{\sin \pi z}+2\pi \sum_{k\geq 1} \left[ \frac{1}{\sin(2k\pi \tau-\pi z)} - \frac{1}{\sin(2k\pi \tau+\pi z)}\right]. $$ 

\section{  Weierstrass's function \ $\wp(z)$\ and infinite products}
The Weierstrass's function \ $\wp(z) = \wp(z;\omega,\omega')$\ is an elliptic function of order two with two primitive periods \ $(2\omega, 2\omega')$ \ verifying\ $\tau = \frac{\omega'}{\omega},\ \Im \tau  > 0$. \\  
Notice that general case of two  periods $(2\omega,2\omega')$ (with $\tau=\frac{\omega'}{\omega}$ ) may be easily deduced by the following homogeneity relations for arbitrary $t \neq 0$ : 
$$\wp(tz,t\omega,t\omega') = t^{-2} \wp(z,\omega,\omega'), \quad \wp'(tz,t\omega,t\omega') = t^{-3} \wp'(z,\omega,\omega').$$
In the sequel we will consider the Weierstrass elliptic function \ $\wp(z) = \wp(z,1,\tau)$\ with two primitive periods $2,2\tau$ one of them is real and the other imaginary, see [6, p.45].
Throughout this paper we will take \ $\omega = 1, \ \omega' = \tau$\ and in order to avoid any ambiguity \ $\wp(z)$\ always denotes \ $\wp(z,\tau)$.

\subsection{ Another expression of \ $\wp(z)$ }
There is a connexion between \ $\wp(z)$\ and theta functions. The following is well known [1] or [4, (4), p.361] for $v = \frac{z}{2\omega} = \frac{z}{2} $:

$$\wp(z) = e_i +\frac{1}{4} \left [ \frac{\theta_{i+1}(v)}{\pi \theta_{i+1}(0)} \frac{ \theta'_{1}(0)}{\theta_{1}(v)}  \right]^2, \quad i=1,2,3.$$
These relations allows us to derive variant infinite products expressing the Weierstrass's function\\

 {\bf Theorem 2-1}   \quad {\it The Weierstrass's function $\wp(z) =\wp(z,\tau)$  with primitive periods $2$ and $2\tau$ verifies the following identities
$$\wp(z) - e_1 = \frac{(\pi \cot \frac{\pi z}{2})^2}{4} \prod_{k\geq 1}  \left[ \frac{ \cot(k\pi \tau- \frac{\pi z}{2})\ \cot(k\pi \tau+ \frac{\pi z}{2})}{\left[\cot(k\pi \tau)\right ]^2}\right]^2 =$$ $$ \frac{(\pi^2 \theta_3(0) \theta_4(0) \cot \frac{\pi z}{2})^2}{4} \prod_{k\geq 1}  \left[{ \cot(k\pi \tau- \frac{\pi z}{2})\ \cot(k\pi \tau+ \frac{\pi z}{2})}\right]^2$$ 
where $e_1 = \wp(1),$ \ and \ $Im z < 2 \ Im \tau$.

We obtain the two other infinite products by permuting the $e_i$} 
$$\wp(z) - e_3 = \frac{\pi^2}{\left(2 \sin \frac{\pi z}{2}\right)^2}\prod_{k\geq 1} \left[ \frac{\sin((k-\frac{1}{2})\pi \tau-\pi \frac{z}{2})}{\sin(k\pi \tau-\pi \frac{z}{2})} \frac{\sin((k-\frac{1}{2})\pi \tau+\pi \frac{z}{2})}{\sin(k\pi \tau+\pi \frac{z}{2})}\right]^2 \left [\frac{\sin(k\pi \tau)}{\sin((k-\frac{1}{2})\pi \tau)}\right]^4,$$
$$\wp(z) - e_2 = \frac{\pi^2}{\left(2 \sin \frac{\pi z}{2}\right)^2}\prod_{k\geq 1} \left[ \frac{\cos((k-\frac{1}{2})\pi \tau-\pi \frac{z}{2})}{\sin(k\pi \tau-\pi \frac{z}{2})} \frac{\cos((k-\frac{1}{2})\pi \tau+\pi \frac{z}{2})}{\sin(k\pi \tau+\pi \frac{z}{2})}\right]^2 \left [\frac{\sin(k\pi \tau)}{\cos((k-\frac{1}{2})\pi \tau)}\right]^4.$$\\

{\bf Proof of Theorem 2-1}\quad Starting from 
$$\wp(z) = e_1 +\frac{1}{4} \left [\pi \frac{\theta_{2}(v)}{ \theta_{2}(0)} \frac{\theta'_{1}(0)}{\theta_{1}(v)}  \right]^2,$$ 
and by [2, Corollary 3-5] which asserts 
 $$\frac{ \theta_1(v,\tau)}{\pi ( \sin \pi v )\ \theta'_1(0,\tau)} =   \left[1 - \left(\frac{\sin \pi v}{\sin k\pi \tau} \right)^2 \right] = \prod_{k\geq 1} \frac{\sin (k\pi \tau-\pi v)\ \sin (k\pi \tau+\pi v)}{\left[\sin (k\pi \tau)\right]^2},$$ 
$$ \frac{\theta_2(v,\tau)}{ (\cos \pi v )\ \theta_2(0,\tau)} =  \left[1 - \left(\frac{\sin \pi v}{\cos k\pi \tau} \right)^2 \right] = \prod_{k\geq 1} \frac{\cos (k\pi \tau-\pi v)\ \cos (k\pi \tau+\pi v)}{\left[\cos (k\pi \tau)\right]^2}.$$
Then, 
$$ \frac{\theta_{2}(v)}{ \theta_{2}(0)} \frac{\pi\theta'_{1}(0)}{\theta_{1}(v)} =  \cot \pi v  \prod_{k\geq 1} \frac{\cos [k\pi \tau+\pi v] \cos [k\pi \tau-\pi v] ( \sin [k\pi \tau])^2}{\sin [k\pi \tau+\pi v] \sin [k\pi \tau-\pi v]( \cos [k\pi \tau])^2}$$
$$=  \cot \pi v  \prod_{k\geq 1} \left[ \frac{\cot(k\pi \tau+\pi v) \cot(k\pi \tau-\pi v)}{(\cot k\pi \tau)^2}\right ] .$$

Moreover,  by the relationship \ $\wp(z) - e_1 = \left(\frac{\sigma_1 z}{\sigma z}\right)^2,$ \  we have for \ $v =\frac{z}{2}$\ (see for example Schwarz [9, p. 8,36])
$$\sigma_1 z = e^{\frac{\eta v^2}{2}}  {( \cos \pi v)} \prod_{k\geq 1}   \frac{ \cos(k\pi \tau-\pi v)\ \cos(k\pi \tau+\pi v)}{\left[\cos (k\pi \tau)\right]^2} $$
 $$\sigma z = e^{\frac{\eta v^2}{2}}  \frac{(2 \sin \pi v)}{\pi} \prod_{k\geq 1}   \frac{ \sin(k\pi \tau-\pi v)\ \sin(k\pi \tau+\pi v)}{\left[\sin (k\pi \tau)\right]^2} .$$
We then derive the  expression
$$\frac{\sigma_1 z}{\sigma z} = \frac{(\pi \cot \pi v)}{2} \prod_{k\geq 1}    \frac{\cot(k\pi \tau-\pi v)\ \cot(k\pi \tau+\pi v)}{\left[\cot (k\pi \tau)\right]^2}, $$
and deduce 
$$\wp(z) - e_1 = \frac{(\pi \cot \frac{\pi z}{2})^2}{4} \prod_{k\geq 1} \left[\frac{ \cot(k\pi \tau-\pi v)\ \cot(k\pi \tau+\pi v)}{\left[\cot (k\pi \tau)\right]^2}\right]^2 =
\frac{(\pi \cot \frac{\pi z}{2})^2}{4} \prod_{k\neq 0}  \left[ \frac{ \cot(k\pi \tau- \frac{\pi z}{2})}{\left[\cot k\pi \tau\right ]}\right]^2.$$

By permutation of the \ $e_i, i=1,2,3$\ we also obtain the other expressions\\ (see Schwarz [9, p.36])
$$\sigma_2 z = e^{\frac{\eta v^2}{2}}   \prod_{k\geq 0}   \frac{ \cos((k-\frac{1}{2})\pi \tau-\pi v)\ \cos((k-\frac{1}{2})\pi \tau+\pi v)}{\left[\cos ((k-\frac{1}{2})\pi \tau)\right]^2} ,$$
$$\sigma_3 z = e^{\frac{\eta v^2}{2}}   \prod_{k\geq 0}   \frac{ \sin((k-\frac{1}{2})\pi \tau-\pi v)\ \sin((k-\frac{1}{2})\pi \tau+\pi v)}{\left[\sin ((k-\frac{1}{2})\pi \tau)\right]^2} ,$$
and then deduce analog infinite products for $$ \wp(z) - e_2 = \left(\frac{\sigma_2 z}{\sigma z}\right)^2 \quad \wp(z) - e_3 = \left(\frac{\sigma_3 z}{\sigma z}\right)^2$$.\\

{\bf Remark 2-2 :}\  {\bf (i)}\quad This expansion as infinite product of elliptic functions may be differently proved. Indeed, starting from the infinite product
$$\sin (\frac{\pi u}{2\omega}) = \frac{\pi u}{2\omega} \prod_{n\geq 0} \left(1 - \frac{ u}{2n \omega}\right) e^{\frac{ u}{2n\omega}} \prod_{n\geq 0} \left(1 + \frac{ u}{2n \omega}\right) e^{-\frac{ u}{2n\omega}}$$H.A. Schwarz [9, p.8] noticed that the sigma function may also be written
$$\sigma z = \frac{2\omega}{\pi} \sin (\pi v) e^{2\eta \omega v^2} \prod_{n\geq 1} \left[\frac{\sin (n\pi \tau-\pi v) \sin (n\pi \tau+\pi v)}{(\sin n\pi \tau)^2}\right]  = e^{2\eta \omega v^2} \frac{2\omega}{\pi} \sin (\pi v) \prod_{n\geq 1} \left(1 - \frac{ (\sin \pi v)^2}{(\sin n\pi \tau)^2}\right),$$
where \ $v = \frac{z}{2\omega}, \ \eta = \frac{\pi^2}{2\omega} \left (\frac{1}{6} + \sum_{n\geq 0} \frac{1}{(\sin n\pi \tau)^2}\right).$\\

{\bf (ii)}\quad Many relations and descriptive properties of \ $\wp(z)$\  may be derived from Theorem 2-1. For example, We may write : 
$$\wp(z+1) = e_1 +\frac{(\pi \tan \pi \frac{z}{2})^2}{4} \prod_{k\neq 0}  \left[ \frac{ \tan(k\pi \tau-\pi \frac{z}{2})}{\cot (k\pi \tau)}\right ]^2.$$
We then find again 
$$\left(\wp(z+1) - e_1\right)\left(\wp(z) - e_1\right) =  \prod_{k\geq 1} \frac{1}{16\left(\cot k\pi \tau \right)^{8}} = (e_1 - e_2)(e_1 - e_3)= \frac{\pi^4}{16} \theta_3^4(0) \theta_4^4(0).$$
By the same way one gets
$$(e_3 - e_2)(e_3 - e_1) = \prod_{k\geq 1} \frac{1}{16\left(\cot \frac{k\pi}{ \tau} \right)^{8}}, \quad  (e_2 - e_1)(e_2 - e_3) = \prod_{k\geq 1} \frac{1}{16\left(\cot \frac{k\pi}{ 1+\tau} \right)^{8}}.$$

\subsection{Logarithmic differentiation of $\wp(z) - e_1$}
Taking the logarithmic differentiation in Theorem 2-1 we found the following\\ 

{\bf Corollary 2-3}   \quad {\it Let \ $\wp'(z)$\ the derivative  of the Weierstrass function relative to the periods \ $(2,2\tau)$. \ We then have}
$$ \frac{\wp'(z)}{\wp(z) - e_1} =  -\frac{2\pi}{\sin \pi z}+2\pi \sum_{k\neq 0} \left[ \frac{1}{\sin(2k\pi \tau-\pi z)}\right] = -\frac{\sigma(2z)}{\sigma^2(z) \sigma^2_1(z)} .$$\\
 
Indeed, that result follows from Theorem 2-1 and the identity $$\frac{d \cot x }{d x}\frac{1}{\cot x} = - \frac{2}{\sin 2x}.$$
Then 
$$\frac{\wp'(z)}{\wp(z) - e_1} =  -\frac{4\pi}{2\sin \pi z} + \frac{4\pi}{2} \sum_{k\geq 1} \left[ \frac{1}{\sin(2k\pi \tau-\pi z)} - \frac{1}{\sin(2k\pi \tau+\pi z)}\right]= \sum_{k} \left[ \frac{2\pi}{\sin(2k\pi \tau-\pi z)}\right].$$ 
We may expressed it otherwise (for example : see Lawden [7, p.161]) 
$$\frac{\wp'(z)}{\wp(z) - e_1} = 2 \left(\frac{\sigma'_1}{\sigma_1}(z) -\frac{\sigma'}{\sigma}(z)\right) = -2 \frac{\sigma_2(z) \sigma_3(z)}{\sigma(z) \sigma_1(z)} = -\frac{\sigma(2z)}{\sigma^2(z) \sigma^2_1(z)} =  2\zeta(z+1)-2\zeta(z)-2\eta.$$ \\ 

Concerning the sigma functions we derive the following results (here \ $\eta = \frac{\pi^2}{2} \left (\frac{1}{6} + \sum_{n\geq 0} \frac{1}{(\sin n\pi \tau)^2}\right)$) using their infinite product expansions [10, T.2, p.246] and their logarithmic derivatives

\bigskip 

{\bf Corollary 2-4}   \quad {\it The elliptic sigma functions and their derivatives relative to periods \ $(2,2\tau)$\ verify the identities}
$$\frac{\sigma'}{\sigma}(z) = \eta z + \frac{\pi}{2} \cot(\frac{\pi z}{2}) + \frac{\pi}{2}\sum_{k\geq 1} \left[\cot(k\pi \tau+\frac{\pi z}{2}) - \cot(k\pi \tau-\frac{\pi z}{2})  \right]=$$
$$\eta z + \frac{\pi}{2} \cot(\frac{\pi z}{2}) + \frac{\pi}{2}\sum_{k\geq 1}  \frac{\sin \pi z}{-\cos^2(\frac{\pi z}{2}) + \cos^2(k\pi \tau)},$$
$$\frac{\sigma'_1}{\sigma_1}(z) = \frac{\sigma'}{\sigma}(z+1) = \eta z -\frac{\pi}{2} \tan(\frac{\pi z}{2})- \frac{\pi}{2}\sum_{k\geq 1} \frac{\sin \pi z}{-\sin^2(\frac{\pi z}{2}) + \cos^2(k\pi \tau)},$$ 
$$\frac{\sigma'_2}{\sigma_2}(z) = \eta z -\frac{\pi}{2} \sum_{k\geq 0} \left[\tan((k-\frac{1}{2})\pi \tau+\frac{\pi z}{2}) - \tan((k-\frac{1}{2})\pi \tau+\frac{\pi z}{2})  \right] = $$
$$\eta z - \frac{\pi}{2}\sum_{k\geq 1} \frac{\sin \pi z}{-\sin^2(\frac{\pi z}{2}) + \cos^2((k-\frac{1}{2})\pi \tau)},$$
$$\frac{\sigma'_3}{\sigma_3}(z) = \frac{\sigma'_2}{\sigma_2}(z+1) = \eta z + \frac{\pi}{2}\sum_{k\geq 1}  \frac{\sin \pi z}{-\cos^2(\frac{\pi z}{2}) + \cos^2((k-\frac{1}{2})\pi \tau)}.$$

\bigskip 

{\bf Corollary 2-5}   \quad {\it The elliptic sigma functions and their derivatives relative to periods \ $(2,2\tau)$\ verify the identities}
$$\frac{\wp'(z)}{2(\wp(z) - e_1)} = \frac{\sigma'_1}{\sigma_1}(z) - \frac{\sigma'}{\sigma}(z) =  \pi \sum_{k\geq 1} \left[ \frac{1}{\sin(2k\pi \tau-\pi z)} -  \frac{1}{\sin(2k\pi \tau+\pi z)}\right],$$
$$ \frac{\sigma'_2}{\sigma_2}(z) - \frac{\sigma'_3}{\sigma_3}(z) = \pi \sum_{k\geq 0} \left[ \frac{1}{\sin((2k-1)\pi \tau-\pi z)}- \frac{1}{\sin((2k-1)\pi \tau+\pi z)}\right].$$
$$\frac{\sigma'_1}{\sigma_1}(z) + \frac{\sigma'}{\sigma}(z) = 2\eta z + {\pi} \cot({\pi z}) + {\pi}\sum_{k\geq 1} \left[\cot(2k\pi \tau+\pi z) - \cot(2k\pi \tau-\pi z)\right],$$
$$ \frac{\sigma'_2}{\sigma_2}(z) + \frac{\sigma'_3}{\sigma_3}(z) =  2\eta z + {\pi}\sum_{k\geq 0} \left[\cot((2k-1)\pi \tau+\pi z) - \cot((2k-1)\pi \tau-\pi z)\right].$$\\

The next result provide other expressions for the Weierstrass function as well as for its derivative
\bigskip 

{\bf Corollary 2-6}   \quad {\it The Weierstrass function relative to periods \ $(2,2\tau)$\ verifies the identities}
$$\wp(z,\tau)-e_2(\tau) = \pi^2 \sum_{k} \left[ \frac{1}{\sin(2k\pi \tau+\pi z)}\right]\ \sum_{k} \left[ \frac{1}{\sin(\frac{-2k\pi}{ \tau}+\pi z)}\right],$$
$$\wp(z,\tau)-e_3(\tau) = \pi^2 \sum_{k} \left[ \frac{1}{\sin(2k\pi \tau+\pi z)}\right]\ \sum_{k} \left[ \frac{1}{\sin(\frac{-2k\pi}{ \tau+1}+\pi z)}\right],$$
$$\wp(z,\tau)-e_1\tau) = \pi^2 \sum_{k} \left[ \frac{1}{\sin(\frac{2k\pi}{ \tau+1}+\pi z)}\right]\ \sum_{k} \left[ \frac{1}{\sin(\frac{-2k\pi}{ \tau}+\pi z)}\right].$$
$$\wp'(z,\tau) = -2\pi^3 \sum_{k} \left[ \frac{1}{\sin(2k\pi \tau+\pi z)}\right]\ \sum_{k} \left[ \frac{1}{\sin(\frac{-2k\pi}{ \tau}+\pi z)}\right] \sum_{k} \left[ \frac{1}{\sin(\frac{-2k\pi}{ \tau+1}+\pi z)}\right].$$\\

Indeed, one has since [4, p.368] \ $e_1(\frac{-1}{\tau}) = e_3(\tau)$ \ and \ $\ e_1(\frac{-1}{\tau+1}) = e_2(\tau)$\ then 
$$\frac{\wp'(z)}{\wp(z) - e_3} = -2\pi \sum_{k} \left[ \frac{1}{\sin(\frac{-2k\pi}{ \tau}+\pi z)}\right],\quad \frac{\wp'(z)}{\wp(z) - e_2} = -2\pi \sum_{k} \left[ \frac{1}{\sin(\frac{-2k\pi}{ \tau+1}+\pi z)}\right].$$
Moreover, since $$\frac{\wp'^2(z)}{(\wp(z) - e_1) (\wp(z) - e_3)} = \left[\frac{\wp'(z)}{(\wp(z) - e_1)}\right] \ \left[\frac{\wp'(z)}{(\wp(z) - e_3)}\right] = 4 \left[\wp(z) - e_2\right]$$ it easily follows the first identity of Corollary 2-6. On the other hand, since \ $e_2(1+\tau) = e_3(\tau), \ e_1(\frac{-1}{\tau+1}) = e_2(\tau)$\ we then obtain the two other identities.\\

On the other hand, by [9, p.13, (1.)]  the Weierstrass function $\wp(z,\tau),$ with primitive periods $2,2\tau$ and the sigma function are related
$$-\frac{\sigma(u+v) \sigma(u-v)}{\sigma^2(u) \sigma^2(v)} = \wp(u)-\wp(v).$$
Since \ $\sigma(u) = \frac{2}{\pi \theta'_1(0,\tau)} e^{\eta u^2/2} \theta_1(\frac{u}{2},\tau)$,\ then
$$\wp(u)-\wp(v) = (\pi \theta'_1(0,\tau))^2 \frac{\theta_1(\frac{u+v}{2},\tau) \theta_1(\frac{u-v}{2},\tau)}{\left(\theta_1(\frac{u}{2},\tau) \theta_1(\frac{v}{2},\tau)\right)^2}.$$
In the limit case when \ $v \rightarrow u$ after dividing both sides by $v-u$, we find
$$\wp'(u) = - \frac{\sigma(2u)}{\sigma^4(u)} = - [\pi \theta'_1(0,\tau)]^3 \frac{\theta_1(2u,\tau)}{[\theta_1(u,\tau)]^4}.$$
However, using again infinite products of theta functions we may express $\wp'(u)$ as infinite product.

\bigskip 
 {\bf Theorem 2-7}  \quad {\it The derivative of Weierstrass's function $\wp'(u) =\wp'(u,\tau)$  with primitive periods $2$ and $2\tau$ verifies the following for \ $v =\frac{u}{2}$}

$$\wp'(u) = -\frac{\sin 2\pi v}{ (\sin\pi v)^4} \prod_{k\geq 1} \frac{\sin(k\pi \tau+2\pi v) \sin(k\pi \tau-2\pi v) (\sin(k\pi \tau))^6}{[\sin(k\pi \tau+\pi v) \sin(k\pi \tau-\pi v)]^4}.$$\\

{\bf Proof of Theorem 2-7}\quad Indeed, by [2, Cor. 3-5] 
$$ \frac{\theta_1(v,\tau)}{ (\pi \sin \pi v )\ \theta'_1(0,\tau)} =  \ \prod_{k\geq 1} \left[1 - \left(\frac{\sin \pi v}{\sin k\pi \tau} \right)^2 \right] = \prod_{k\geq 1} \frac{\cos 2\pi v - \cos 2k\pi \tau}{1 - \cos 2k\pi \tau},$$
$$[ \frac{\theta_1(v,\tau)}{ (\pi \sin \pi v )\ \theta'_1(0,\tau)} ]^4 =  \ \prod_{k\geq 1} \left[1 - \left(\frac{\sin \pi v}{\sin k\pi \tau} \right)^2 \right]^4 = \prod_{k\geq 1} \left [\frac{\cos 2\pi v - \cos 2k\pi \tau}{1 - \cos 2k\pi \tau}\right]^4.$$
Thus
$$ \frac{\theta_1(2v,\tau)}{ \theta'_1(0,\tau)}\frac{(\theta'_1(0,\tau))^4}{[\theta_1(v,\tau)]^4} = - \pi^3 \ \wp'(z) = \frac{\pi^4 \sin 2\pi v }{\pi [\sin \pi v]^4} \prod_{k\geq 1} {\frac{\frac{\cos 4\pi v - \cos 2k\pi \tau}{1 - \cos 2k\pi \tau}}{\left[\frac{\cos 2\pi v - \cos 2k\pi \tau}{1 - \cos 2k\pi \tau}\right]^4}}$$
$$ = \frac{ \pi^3 \sin 2\pi v }{ (\sin \pi v)^4} \prod_{k\geq 1} \frac{[\cos 4\pi v - \cos 2k\pi \tau] [1 - \cos 2k\pi \tau]^3}{[\cos 2\pi v - \cos 2k\pi \tau]^4}$$
$$ = \frac{\pi^3 \sin 2\pi v}{ (\sin\pi v)^4} \prod_{k\geq 1} \frac{\sin(k\pi \tau+2\pi v) \sin(k\pi \tau-2\pi v) (\sin(k\pi \tau))^6}{[\sin(k\pi \tau+\pi v) \sin(k\pi \tau-\pi v)]^4}. $$

\section{The functions \ $\xi(u,\tau)$}

Following Tannery [10, T.2, p. 168] it is often appropriate to introduce the quotients that can be formed by means of two functions $\sigma (z,\tau)$ relating to the same variable and the same periods. Thus the four functions \ $\sigma, \sigma_1, \sigma_2, \sigma_3$\  generate twelve functions. Indexes \ $\alpha, \beta,\gamma$ \ are selected according to the convention : they are different and be chosen between the values \ $0, 1, 2, 3$\ 
$$\xi_{\alpha 0}(u) = \frac{\sigma_\alpha}{\sigma}(u) = \sqrt{\wp(u) - e_\alpha},\qquad 
\xi_{0 \alpha}(u) = \frac{\sigma}{\sigma_\alpha}(u) = \frac{1}{\sqrt{\wp(u) - e_\alpha}},$$ $$ 
\xi_{\beta \gamma}(u) = \frac{\sigma_\beta}{\sigma_\gamma}(u) = \frac{\sqrt{\wp(u) - e_\beta}}{\sqrt{\wp(u) - e_\gamma}}.$$
We will simply write \ $\xi_{\beta \gamma}(u) $\ instead \ $\xi_{\beta \gamma}(u,\tau)$\ when there is no ambiguity.
These functions are even or odd depending on whether they contain the index 0 or not. They are algebraic function of \ $\wp(u)$\ and have a single pole as simple singularity. Moreover, these functions are doubly periodic and verify many algebraic relations see [10, T.2, p. 168-187].
Observe that we may deduce from the above :
$$\wp'(u) = - 2 \frac{\sigma_1}{\sigma}(u)\ \frac{\sigma_2}{\sigma}(u)\ \frac{\sigma_3}{\sigma}(u) = -2 \xi_{\alpha 0}(u)\ \xi_{\beta 0}(u)\ \xi_{ \gamma 0}(u).$$
These relations yield in particular, [10, T.2, p. 171]
$$\xi'_{\alpha 0}(u) = \frac{\wp'(u)}{2\sqrt{\wp(u) - e_\alpha}} = - {\sqrt{\wp(u) - e_\beta}}\ {\sqrt{\wp(u) - e_\gamma}} = - \xi_{\beta 0}(u)\  \xi_{\gamma 0}(u),$$
$$\xi'_{0 \alpha}(u) = \xi_{\beta \alpha}(u)\ \xi_{\gamma \alpha}(u),\quad \xi'_{\beta \gamma}(u) =  - (e_\beta - e_\gamma) \xi_{0 \gamma}(u)\  \xi_{\alpha \gamma}(u), $$ $$ \xi_{\beta 0}(u) \ \xi_{\gamma 1}(u) = \frac{-\wp'(u)}{2(\wp(u)- e_1)}.$$
Moreover, \ $\xi_{\alpha 0}$\ is solution of the differential equation
$$\left(\frac{dy}{du}\right)^2 = (e_\alpha - e_\beta + y^2)\ (e_\alpha - e_\gamma + y^2).$$
The derivatives of the Weierstrass functions yield
$$3 \wp(u) = \xi^2_{\alpha 0}(u) + \xi^2_{\beta 0}(u) + \xi^2_{\gamma 0}(u), \quad \frac{\wp''(u)}{\wp'(u)} = 2\xi'_{\alpha 0}(u) + 2\xi'_{\beta 0}(u) + 2\xi'_{\gamma 0}(u).$$
Notice for \ $ \omega_1 = 1, \omega_2 = 1+\tau, \omega_3 = \tau$\ we derive
$$\xi_{\alpha 0}(\omega_\beta) = \sqrt{e_\beta-e_\alpha}, \quad \xi_{\beta \gamma}(\omega_\alpha) = \frac{\sqrt{e_\alpha-e_\beta}}{\sqrt{e_\alpha-e_\gamma}}.$$
The periods are
$$k = \xi_{2 1}(\omega_3), \qquad k' = \xi_{2 3}(\omega_1).$$

Many other interesting properties concerning the functions \ $\xi_{\beta \gamma}(u)$\ (homogeneity, variations, growing,...) can be found in [10, T.2, p. 170-190]. In [10] we will also observe a rigorous and systematic study of these functions as well as their periodicity, their relationships with other elliptic functions and connections with theta functions of Jacobi. The following is particularly interesting\\

{\bf Theorem 3-1}   \quad {\it The function \ $\xi_{\alpha 0}(u,\tau) =\xi_{\alpha 0}(u)$\ verify the identities}  
$$\frac{\xi'_{1 0}(u)}{\xi_{1 0}(u)} = \sqrt{\left(\frac{e_1 - e_2}{\xi^2_{1 0}(u)} + 1\right)\ \left(\frac{e_1 - e_3}{\xi^2_{1 0}(u)} + 1\right)} = -\pi \sum_{k\neq 0} \left[ \frac{1}{\sin(2k\pi \tau+\pi z)}\right],$$
$$\frac{\xi'_{2 0}(u)}{\xi_{2 0}(u)} = \sqrt{\left(\frac{e_2 - e_1}{\xi^2_{2 0}(u)} + 1\right)\ \left(\frac{e_2 - e_3}{\xi^2_{2 0}(u)} + 1\right)} = -\pi \sum_{k\neq 0} \left[ \frac{1}{\sin(\frac{-2k\pi}{\tau+1}+\pi z)}\right],$$
$$\frac{\xi'_{3 0}(u)}{\xi_{3 0}(u)} = \sqrt{\left(\frac{e_3 - e_2}{\xi^2_{3 0}(u)} + 1\right)\ \left(\frac{e_3 - e_1}{\xi^2_{3 0}(u)} + 1\right)} = -\pi \sum_{k\neq 0} \left[ \frac{1}{\sin(\frac{-2k\pi}{\tau}+\pi z)}\right].$$\\

{\bf Proof of Theorem 3-1}\quad This theorem is a direct consequence of Corollaries 2-3 and 2-5 .

\section{Transformations of order \ $n$\  }
Then the $n$-transformation theory of functions deals with the relations
between functions belonging to different pairs of primitive periods : \ $(2,2\tau),\quad (2,2n\tau).$\\

If we suppose \ $n$\  is not a prime number it will be the product of two or more odd primes, and the transformation will break up into distinct transformations each of which may be separately considered. We therefore now assume \ $n$\ an odd prime: the modular equation is in this case an irreducible equation of the order \ $n+1$.\\ 
Then, it is convenient to restrict the attention to the case \ $n$\ an odd number. Thus, we are going to study the transformations whose order is odd
and positive, transformations which are reduced to those that we change the period \ $(2,\tau)$\  into \ $(\frac{2}{n},\tau)$\  without changing $\tau$.\ This corresponds to changing \ $v$\ to \ $nv$\ and \ $\tau$\ to \ $n\tau$.\\ 
 We shall always assume
$$Im(\tau) > 0,$$
 Observe by Lawden [7, p.252], Enneper [4, p.240] or Roy [9, p.89] that any transformation of order $ n > 1$  may be represented as a product of transformations of first order and of transformations of higher order with matrix 
$$ M = \left(
\begin{array}{clrr} 
1 & 0 \\
0 & n 
\end{array}\right)
$$ 
Moreover, any transformation \ $\tau' = n\tau$\ can be separated into a product when $n$ has prime factors. Therefore, we only study the case of transformation when $n$ is a prime and limit our study for this type of matrix. 

\subsection{Transformation of \ $\wp(z,\tau)$}
Theorem 2-1 provided an expansion of the Weierstrass function \ $\wp(z,\tau)$ \ relative to periods \ $2,2\tau)$\ as infinite product 
$$\wp(z) = e_1 +\frac{(\pi \cot \frac{\pi z}{2})^2}{4} \prod_{k\neq 0}  \left[ \frac{ \cot(k\pi \tau- \frac{\pi z}{2})}{\cot k\pi \tau}\right ]^2,$$ 
where $e_1 = \wp(1)$. It allows to derive in particular a $n$-order transformation formula for $i = 1,2,3$.  

That permits to deduce various identities as
\begin{equation} \wp(nz,n\tau) - e_j(n\tau) = [\wp(z,\tau) - e_j(\tau)] \prod_{m=1}^{\frac{n-1}{2}}\left[ \frac{\wp(z) - \wp(\frac{m}{n}+\omega_j)}{\wp(z) - \wp(\frac{m}{n})} \right]^2\end{equation}
where \ $j=1,2,3, \ \omega_1=1,\ \omega_2=1+\tau,\ \omega_3=\tau,$ \ this means \\ $e_1 = \wp(1), \ e_2 = \wp(1+\tau),\ e_3 = \wp(\tau).$\\

More precisely, one gets the following connection between $\wp(nz,n\tau)$\ and \ $\wp(z,\tau)$ \ derived from Theorem 2-1\\

 {\bf Theorem 4-1}   \quad {\it Let $n$ be an odd integer and consider the Weierstrass's function \ $\wp(nz,n\tau)$\  with primitive periods $2$ and $2n\tau$, then the following identity holds
$$\wp(nz,n\tau) - e_1(n\tau) = \left(\frac{4}{\pi^2}\right)^{n-1} \prod_{k\geq 1} \left[\frac{(\cot k\pi \tau)^n}{\cot (k n\pi \tau)} \right]^{4}\ \prod^{n-1}_{m=0}  \left[\wp(z+\frac{m}{n},\tau) - e_1(\tau)\right] =$$ 
$$\left(\frac{4}{\pi^4}\right)^{n-1} \frac{\theta_3^2(0,n\tau) \theta_4^2(0,n\tau)}{\left[\theta_3^2(0,\tau) \theta_4^2(0,\tau)\right]^n} \ \prod^{n-1}_{m=0}  \left[\wp(z+\frac{m}{n},\tau) - e_1(\tau)\right],$$ 
$$\wp(nz,n\tau) - e_2(n\tau) = \left(\frac{4}{\pi^4}\right)^{n-1} \frac{\theta_2^2(0,n\tau) \theta_4^2(0,n\tau)}{\left[\theta_2^2(0,\tau) \theta_4^2(0,\tau)\right]^n} \ \prod^{n-1}_{m=0}  \left[\wp(z+\frac{m}{n},\tau) - e_2(\tau)\right],$$ 
$$\wp(nz,n\tau) - e_3(n\tau) = \left(\frac{4}{\pi^4}\right)^{n-1} \frac{\theta_2^2(0,n\tau) \theta_3^2(0,n\tau)}{\left[\theta_2^2(0,\tau) \theta_3^2(0,\tau)\right]^n} \ \prod^{n-1}_{m=0}  \left[\wp(z+\frac{m}{n},\tau) - e_3(\tau)\right],$$ 
where \ $e_j(\tau)$ \ are the zeros of \ $ \wp'(z,\tau),\ i=1,2,3$\  and } \ $Im z < 2\ Im \tau$. \\

{\bf Proof of Theorem 4-1}\quad Prove at first for $i=1$. Notice by Theorem 2-1
$$\wp(nz,n\tau) = e_1(n\tau) +\frac{(\pi \cot \pi nv)^2}{4} \prod_{k\geq 1} \left[ \frac{\cot(kn\pi \tau+n\pi v) \cot(kn\pi \tau-n\pi v)}{(\cot kn\pi \tau)^2}\right ]^2.$$
We Start from the classical trigonometric product formulas valid for $n$ odd integer
$$\sin(nz) = 2^{n-1} \sin(z) \sin(z+\frac{\pi}{n})  \sin(z+\frac{2\pi}{n}).... \sin(z+\frac{(n-1)\pi}{n}) = 2^{n-1}\ \prod^{n-1}_{m=0} \sin(z+\frac{m\pi}{n}).$$
$$\cos(nz) = (-1)^{\frac{n-1}{2}}\ 2^{n-1} \prod^{n-1}_{m=0} \cos(z+\frac{m\pi}{n}),\quad 
\cot(nz) = (-1)^{\frac{n-1}{2}}\  \prod^{n-1}_{m=0} \cot(z+\frac{m\pi}{n}).$$
Thus we derive the expression
$$\wp(nz,n\tau) = e_1(n\tau) + \prod^{n-1}_{m=0}\frac{(\pi \cot (\pi v+\frac{m\pi}{n}))^2}{4} \prod^{n-1}_{m=0} \prod_{k\geq 1} \left[ \frac{\cot(k\pi \tau+\pi v+\frac{m\pi}{n}) \cot(k\pi \tau- \pi v+\frac{m\pi}{n})}{(\cot k  n\pi \tau)^2}\right ]^2.$$
However, we may write for any integer $m \in [0,n-1]$
$$\wp(z+\frac{m}{n},\tau) - e_1(\tau) = \frac{(\pi \cot (\pi v+\frac{m\pi}{n}))^2}{4} \prod_{k\geq 1} \left[ \frac{\cot(k\pi \tau+\pi v+\frac{m\pi}{n}) \cot(k\pi \tau- \pi v+\frac{m\pi}{n})}{(\cot (k \pi \tau)^2}\right ]^2,$$ 
which means
$$\left(\frac{4}{\pi^2}\right)^n \prod_{k\geq 1} \left[\frac{(\cot k\pi \tau)^n}{\cot (k n\pi \tau)} \right]^{4}\ \prod^{n-1}_{m=0}  \left[\wp(z+\frac{m}{n},\tau) - e_1(\tau)\right] = \frac{4}{\pi^2} \left[\wp(nz,n\tau) - e_1(n\tau)\right].$$ 
Turn now to the other cases $j=2,3$, by the same way we deduce analog decomposition of \ $\wp(nz,n\tau) - e_3(n\tau),\ \wp(nz,n\tau) - e_2(n\tau)$\  from Theorem 2-1 since
$$\wp(z) = e_3 + \frac{\pi^2}{\left(2\sin \frac{\pi z}{2}\right)^2}\prod_{k\geq 1} \left[ \frac{\sin((k-\frac{1}{2})\pi \tau-\pi \frac{z}{2})}{\sin(k\pi \tau-\pi \frac{z}{2})}\right]^2 \left [\frac{\sin(k\pi \tau)}{\sin((k-\frac{1}{2})\pi \tau)}\right]^2,$$
$$\wp(z) = e_2 + \frac{\pi^2}{\left(2\sin \frac{\pi z}{2}\right)^2}\prod_{k\geq 1} \left[ \frac{\cos((k-\frac{1}{2})\pi \tau-\pi \frac{z}{2})}{\sin(k\pi \tau-\pi \frac{z}{2})}\right]^2 \left [\frac{\sin(k\pi \tau)}{\cos((k-\frac{1}{2})\pi \tau)}\right]^2.$$\\

Notice under the action of the modular group $\Gamma_0$ the permutation between the $e_j$ does not change the Weierstrass function $\wp(z,\tau).$ Indeed by changing $\tau$ into $\tau+1$ or into $-\frac{1}{\tau}$ it yields [4, p. 365]: $$e_1(\frac{-1}{\tau}) = e_3(\tau), \ e_1(\frac{-1}{1+\tau}) = e_2(\tau),\ e_1(\frac{-1}{\tau}) = e_3(\tau), \ e_3(1+\tau) = e_2(\tau).$$
Therefore, we find analog expansions for \ $ \wp(nz,n\tau) - e_2(n\tau)$\ and \ $\wp(nz,n\tau) - e_3(n\tau).$ \\

Since  $$\wp'^2(z) = 4 (\wp(z)-e_1)  (\wp(z)-e_2)  (\wp(z)-e_3),\quad and \quad {\theta'}_1(0) = \pi {\theta}_2(0) {\theta}_3(0) {\theta}_4(0)$$ then by Remark 2-2 (ii) we may deduce the following $n$-decomposition of \ $\wp'(nz,n\tau)$ .

\bigskip 

 {\bf Corollary 4-2}   \quad {\it Let $n$ be an odd integer, then the following identity holds}
$$\wp'(nz,n\tau) = \left(\frac{4}{\pi^4}\right)^{n-1}\frac{{\theta'}_1^2(0,n\tau)}{\left[{\theta'}_1^2(0,\tau)\right]^n} \ \prod^{n-1}_{m=0}  \wp'(z+\frac{m}{n},\tau).$$\\

Among others interesting equalities, we may derive the following identities 
$$\frac{\wp(nz,n\tau) - e_1(n\tau) }{ [\wp(z,\tau) - e_1(\tau)]} = \left(\frac{4}{\pi^4}\right)^{n-1} \frac{\theta_3^2(0,n\tau) \theta_4^2(0,n\tau)}{\left[\theta_3^2(0,\tau) \theta_4^2(0,\tau)\right]^n} \ \prod^{n-1}_{m=1}  \left[\wp(z+\frac{m}{n},\tau) - e_1(\tau)\right] = $$
$$\prod_{m=1}^{n-1}\left[ \frac{\wp(z) - \wp(\frac{m}{n}+1)}{\wp(z) - \wp(\frac{m}{n})} \right]^2 = \left(\frac{4}{\pi^2}\right)^{n-1} \prod_{k\geq 1} \left[\frac{(\cot k\pi \tau)^n}{\cot (k n\pi \tau)} \right]^{4}\ \prod^{n-1}_{m=1}  \left[\wp(z+\frac{m}{n},\tau) - e_1(\tau)\right]. $$

From Theorem 4-1 we also get\\ 
 
 {\bf Corollary 4-3}   \quad {\it Let \ $\wp(z,\tau)$\ be the Weierstrass function and  \ $\wp'(z,\tau)$\ its derivative relative to the periods $(2,2\tau),$\ then the following identities hold for any odd $n$ and $j=1,2,3$}
$$(i) \quad \frac{\wp(nz,n\tau) - e_j(n\tau)}{\wp(z,\tau) - e_j(\tau)} = \prod^{n-1}_{m=1} \left[\frac{\wp(z+\frac{2m}{n}) - e_j(\tau)}{\wp(\frac{2m}{n}) - e_j(\tau)}\right] = \prod^{n-1}_{m=1} [\wp(z+\frac{2m}{n}) - e_j(\tau)] \left[\frac{\sigma}{\sigma_j}(\frac{2m}{n},\tau)\right]^2,$$
$$(ii)  \quad \wp'(nz,n\tau) = 2^{1-n} \wp'(z,\tau) \prod^{n-1}_{m=1} \frac{\wp'(z+\frac{2m\pi}{n})}{\wp'(\frac{2m\pi}{n})},$$
$$(iii)\quad \frac{\wp'(nz,n\tau)}{\wp(nz,n\tau) - e_j(n\tau)} = 2^{1-n} \prod^{n-1}_{m=1} \frac{\wp(\frac{2m}{n},\tau) - e_j(\tau)}{\wp'(\frac{2m}{n},\tau)}\ \prod^{n-1}_{m=0} \frac{\wp'(z+\frac{2m}{n},\tau)}{\wp(z+\frac{2m}{n},\tau) - e_j(\tau)} .$$\\

Indeed, (i) may be deduced from Theorem 4-1 and also by the $n$-order transformation of the quotient of sigma functions (see also [10, T.2, p.215], or [11, ex 9 p.456])
$$\frac{\sigma_j}{\sigma}(nz,n\tau) = \prod^{n-1}_{m=0} \frac{\sigma_j}{\sigma}(z+\frac{2m\pi}{n},\tau) \prod^{n-1}_{m=0}  \frac{\sigma}{\sigma_j}(\frac{2m\pi}{n},\tau),\quad j=1,2,3.$$ 
(ii) is deduced from (i) since $$\wp'(u) = - 2 \frac{\sigma_1}{\sigma}(u)\ \frac{\sigma_2}{\sigma}(u)\ \frac{\sigma_3}{\sigma}(u).$$\\
Moreover, since $n$ is odd then the quotient (ii) over (i) yields (iii).\\
 
 {\bf Corollary 4-4}   \quad {\it Let \ $\wp(z,\tau)$\ be the Weierstrass function and  \ $\wp'(z,\tau)$\ its derivative relative to the periods $(2,2\tau),$\ then the following identities hold for any odd $n$ and $j=1,2,3$
$$(i) \ \frac{\wp'(nz,n\tau)}{\wp(nz,n\tau) - e_j(n\tau)} = \frac{\theta_{j+1}^2(0,n\tau)}{\left[\theta_{j+1}^2(0,\tau)\right]^n}\ \prod^{n-1}_{m=0} \frac{\wp'(z+\frac{2m}{n},\tau)}{\wp(z+\frac{2m}{n},\tau) - e_j(\tau)}.$$
Moreover, one gets for $j=1$}
$$(ii)  \ \sum_{k\neq 0} \left[ \frac{1}{\sin(2kn\pi \tau+n\pi z)}\right] = \frac{\theta_{2}^2(0,n\tau)}{\left[\theta_{2}^2(0,\tau)\right]^n}\ \prod^{n-1}_{m=0} \sum_{k\neq 0} \left[ \frac{1}{\sin(2k\pi \tau+\pi z+ \frac{m\pi}{n})}\right].$$\\

Indeed, (i) is derived from Theorem 4-1 and Corollary 4-2 and using the quotient \ $ \frac{\wp'(nz,n\tau)}{\wp(nz,n\tau)-e_1(n\tau)}.$\ Corollary 2-3 implies (ii). 

\subsection {The $n$-transformations of the functions $\xi(u)$}
We start from the $n$-transformations of the sigma functions. In [10,T.1, p.234] Tannery-Molk proved 
$$ \sigma(nu,n\tau) = e^{-nu\sum_{1}^{n-1} \frac{m}{n}\eta_1+nu\wp(\frac{m}{n})}\ \prod_{0\leq m\leq n-1} \frac{\sigma(u+\frac{m}{n},\tau)}{\sigma(\frac{m}{n},\tau)},$$
$$ \sigma_j(nu,n\tau) = e^{-nu\sum_{1}^{n-1} \frac{m}{n}\eta_1+nu\wp(\frac{m}{n})}\ \prod_{0\leq m\leq n-1} \frac{\sigma_j(u+\frac{m}{n},\tau)}{\sigma_j(\frac{m}{n},\tau)},$$ for $j=1,2,3$ we obtain the quotient
$$\frac{\sigma_j(nu,n\tau)}{\sigma(nu,n\tau)} = \prod_{0\leq m\leq n-1}[ \frac{\sigma(\frac{m}{n},\tau)}{\sigma_j(\frac{m}{n},\tau)}] \ \prod_{0\leq m\leq n-1} \frac{\sigma_j(u+\frac{m}{n},\tau)}{\sigma(u+\frac{m}{n},\tau)}.$$

We then deduce [10, T.2, p.215]
$$\xi_{\alpha 0}(nu,n\tau) = \xi_{\alpha 0}(u)  \prod_{m=1}^{n-1} \frac{\xi_{\alpha 0}(u+\frac{2m}{n}) }{\xi_{\alpha 0}(\frac{2m}{n}) },$$
$$\xi_{\beta \gamma}(nu,n\tau) = \xi_{\beta \gamma}(u)  \prod_{m=1}^{n-1} \frac{\xi_{\beta \gamma}(u+\frac{2m}{n}) }{\xi_{\beta \gamma}(\frac{2m}{n}) }.$$
This shows that $\xi_{\alpha 0}(nu,n\tau))$ is a rational function of $\xi_{\alpha 0}(u)$, as well as $\xi_{\beta \gamma}(nu,n\tau)$ is a rational function of $\xi_{\beta \gamma}(u,\tau)$.\\

We then derive the identity 
$$\frac{\xi'_{\alpha 0}(nu,n\tau)}{\xi_{\alpha 0}(nu,n\tau)} = \frac{\xi'_{\alpha 0}(u,\tau)}{\xi_{\alpha 0}(u,\tau)} + \sum_{m=1}^{n-1} \frac{\xi'_{\alpha 0}(u+\frac{2m}{n},\tau)}{\xi_{\alpha 0}(u+\frac{2m}{n},\tau)} = \sum_{m=0}^{n-1} \frac{\wp'(u+\frac{2m}{n})}{2(\wp(u+\frac{2m}{n}) - e_\alpha)},$$
$$\frac{\xi'_{\beta \gamma}(nu,n\tau)}{\xi_{\beta \gamma}(nu,n\tau)} = \frac{\xi'_{\beta \gamma}(u,\tau)}{\xi_{\beta \gamma}(u,\tau)} + \sum_{m=1}^{n-1} \frac{\xi'_{\beta \gamma}(u+\frac{2m}{n},\tau)}{\xi_{\beta \gamma}(u+\frac{2m}{n},\tau)} = \sum_{m=0}^{n-1} \frac{\wp'(u+\frac{2m}{n})}{2(\wp(u+\frac{2m}{n}) - e_\beta)} - \frac{\wp'(u+\frac{2m}{n})}{2(\wp(u+\frac{2m}{n}) - e_\gamma)}.$$

We obtain the modular relations of the periods
$$l = \xi_{2 1}(\tau, n\tau) = \xi_{2 1}(\omega_3) \prod_{m=1}^{n-1} \frac{\xi_{2 1}(\tau+\frac{2m}{n}) }{\xi_{2 1}(\frac{2m}{n}) } = k^n \prod_{m=1}^{n-1} {\xi^2_{1 2}(\frac{2m}{n}) },$$ $$l' = k'^n \prod_{m=1}^{n-1} \frac{1}{\xi^2_{2 1}(\frac{2m}{n})} = k'^n \prod_{m=1}^{n-1} {\xi^2_{3 2}(\frac{2m}{n})}.$$
The zeros become
$$\frac{\sqrt{e_1(\tau)-e_3(\tau)}}{\sqrt{e_1(n\tau)-e_3(n\tau)}} = \prod_{m=1}^{n-1} \frac{\xi^2_{0 3}(\frac{2m-1}{n})}{\xi^2_{0 3}(\frac{2m}{n})}.$$

By the same way if we replace $\tau$ by $\frac{\tau}{n}$ one obtains the formulas
$$\xi_{\alpha 0}(\frac{u}{n},\frac{\tau}{n}) = \xi_{\alpha 0}(u)  \prod_{m=1}^{n-1} \frac{\xi_{\alpha 0}(u+\frac{2m\tau}{n}) }{\xi_{\alpha 0}(\frac{2m\tau}{n}) },$$
$$\xi_{\beta \gamma}(\frac{u}{n},\frac{\tau}{n}) = \xi_{\beta \gamma}(u)  \prod_{m=1}^{n-1} \frac{\xi_{\beta \gamma}(u+\frac{2m\tau}{n}) }{\xi_{\beta \gamma}(\frac{2m\tau}{n}) }.$$
More generally ([10, T.2, p. 217]), if we replace $\tau$ by $\frac{\tau+2p}{n}$ for any odd $n$ and integer $p$ one gets
$$\xi_{\alpha 0}(\frac{u}{n},\frac{\tau+2p}{n}) = \xi_{\alpha 0}(u)  \prod_{m=1}^{n-1} \frac{\xi_{\alpha 0}(u+\frac{2m(\tau+2p)}{n}) }{\xi_{\alpha 0}(\frac{2m(\tau+2p)}{n}) },$$
$$\xi_{\beta \gamma}(\frac{u}{n},\frac{\tau+2p}{n}) = \xi_{\beta \gamma}(u)  \prod_{m=1}^{n-1} \frac{\xi_{\beta \gamma}(u+\frac{2m(\tau+2p)}{n}) }{\xi_{\beta \gamma}(\frac{2m(\tau+2p)}{n}) }.$$\\

Other applications may be derived from the decomposition of \ $\frac{\wp'(nz,n\tau)}{\wp(nz,n\tau)-e_1(n\tau)}$\ as the following modular identities 

\bigskip 
 {\bf Corollary 4-5}   \quad {\it For \ $0< \Im(\tau) < z $\ and for any odd integer \ $n$\ the following identity holds}
$$\frac{\xi'_{1 0}(nz,n\tau)}{\xi_{1 0}(nz,n\tau)} = \frac{\wp'(nz,n\tau)}{2[\wp(nz,n\tau)-e_1(n\tau)]} = \sum_{m=0}^{n-1} \frac{\wp'(z+\frac{2m}{n})}{2(\wp(z+\frac{2m}{n}) - e_1)}=$$ $$-\pi \sum_{k\neq 0} \left[ \frac{1}{\sin(2kn\pi \tau+n\pi z)}\right] = -\pi \sum_{m=0}^{n-1} \sum_{k\neq 0} \left[ \frac{1}{\sin(2k\pi \tau+\pi (z+\frac{2m}{n}))}\right].$$ \\

To prove that, take the logarithmic differentiation of \ $\xi_{1 0}(nz,n\tau)$\ see [10, T.2, p.215] :
 $$\xi_{1 0}(nz,n\tau) = \xi_{1 0}(z,\tau) \prod_{m=1}^{n-1} \frac{\xi_{1 0}(z+\frac{2m}{n},\tau)}{\xi_{1 0}(\frac{2m}{n},\tau)},$$
which yields
$$\frac{\xi'_{1 0}(nz,n\tau)}{\xi_{1 0}(nz,n\tau)} = \frac{\xi'_{1 0}(z,\tau)}{\xi_{1 0}(z,\tau)} \sum_{m=1}^{n-1} \frac{\xi'_{1 0}(z+\frac{2m}{n},\tau)}{\xi_{1 0}(z+\frac{2m}{n},\tau)}.$$
Therefore, Corollary 2-3 implies Corollary 4-4.\\  

 {\bf Corollary 4-6}   \quad {\it For \ $0< \Im(\tau) < z $\ and for any odd integer \ $n$\ the following identities hold}
$$\sum_{k\neq 0} \left[ \frac{1}{\sin(2kn\pi \tau+n\pi z)}\right] = \left(4\pi\right)^{n-1} \left[ \prod_{k\geq 1} \frac{ (\cot \frac{k\pi}{\tau})^{8n}}{  \cot (\frac{k n\pi}{\tau})^8}\right] \prod_{m=0}^{n-1} \left[ \sum_{k\neq 0} \frac{1}{\sin(2k\pi \tau+\pi (z+\frac{m}{n}))}\right] =$$ $$2^{1-n}\sum_{k\neq 0} \prod_{m=1}^{n-1} \left[ \frac{1}{\sin(2k\pi \tau+\pi (z+\frac{2m}{n}))}\right] =  \frac{\theta_{2}^2(0,n\tau)}{\left[\theta_{2}^2(0,\tau)\right]^n}\ \prod^{n-1}_{m=0} \sum_{k\neq 0} \left[ \frac{1}{\sin(2k\pi \tau+\pi z+ \frac{m\pi}{n})}\right].$$\\

To prove Corollary 4-6, recall by Corollary 2-3 we have  
$$\frac{\wp'(z)}{\wp(z) - e_1} =  -2\pi \sum_{k\neq 0} \left[ \frac{1}{\sin(2k\pi \tau+\pi z)}\right] .$$
Since \ $n$\ is odd, by Theorem 2-1, Corollaries 2-3 and 4-3 we get

$$\frac{\wp'(nz,n\tau)}{\wp(nz,\tau) - e_1(n\tau)} =  (-2\pi) \sum_{k\neq 0} \left[ \frac{1}{\sin(2kn\pi \tau+n\pi z)}\right] = \frac{\theta_{j+1}^2(0,n\tau)}{\left[\theta_{j+1}^2(0,\tau)\right]^n}\ \prod^{n-1}_{m=0} \frac{\wp'(z+\frac{2m}{n},\tau)}{\wp(z+\frac{2m}{n},\tau) - e_j(\tau)} = $$ $$ 2^{1-n}\prod^{n-1}_{m=0} \frac{\wp(\frac{m}{n},\tau) - e_j(\tau)}{\wp'(\frac{m}{n},\tau)} \prod_{m=0}^{n-1} (-2\pi)^{n} \sum_{k\neq 0} \left[ \frac{1}{\sin(2k\pi \tau+\pi (z+\frac{m}{n}))}\right]. $$
Thus, the left equality is proved. To prove the right equality, it suffices to use the decomposition of the sinus function  $$\sin(2kn\pi \tau+n\pi z) = 2^{n-1} \prod_{m=0}^{n-1} \sin(2k\pi \tau+\pi (z+\frac{m}{n})).$$
We then obtain an equality allowing us to invert the sum and the product : 
$$2^{1-n}\sum_{k\neq 0} \prod_{m=1}^{n-1} \left[ \frac{1}{\sin(2k\pi \tau+\pi (z+\frac{2m}{n}))}\right] =  \frac{\theta_{2}^2(0,n\tau)}{\left[\theta_{2}^2(0,\tau)\right]^n}\ \prod^{n-1}_{m=0} \sum_{k\neq 0} \left[ \frac{1}{\sin(2k\pi \tau+\pi z+ \frac{m\pi}{n})}\right].$$\\

\vspace{0.5cm}
{\bf References}\\

[1] \ P. Appell, E. Lacour, \quad {\it Fonctions elliptiques et applications}\quad Gauthiers-Villard ed., Paris (1922).\\

[2] \ A.R. Chouikha \quad {\it Functions related to Jacobi Theta Functions and applications}\quad https://hal.archives-ouvertes.fr/hal-03170818. (2021)\\

[3] \ A. Enneper, Elliptische Functionen: Theorie und Geschichte, Louis Nebert, Halle, 1890.\\

[4] \ A.Erdelyi, W.Magnus, F.Oberhettinger,  F.Tricomi, \quad {\it  Higher transcendental functions}\quad  Vol II. Based on notes left by H. Bateman. Robert E. Krieger Publish. Co., Inc., Melbourne, Fla., (1981).\\

[5] \ L. Kiepert, \quad {\it Ueber Theilung und Transformation der elliptischen Functionen.}\quad  Mathematische Annalen t. XXVI, 1885.\\

[6] \ S. Lang, \quad {\it Elliptic functions}\quad Addison-Wesley, Springer, 1970.\\ 

[7] \ D.F. Lawden, \quad {\it Elliptic functions and applications}, Springer-Verlag, vol 80, 1980.\\ 

[8] \ R. Roy,\quad {\it Elliptic and Modular Functions from Gauss to Dedekind to Hecke}, Cambridge Univ. Press, 2017.\\

[9] \ H.A. Schwarz, \quad {\it Formeln und Lehrsatze zum Gebrauche der elliptischen functionen}. Berlin, 1893,\\ 

[10] \ J. Tannery et J. Molk,\quad {\it Elements de la theorie des Fonctions Elliptiques}, tome 2, Gauthier-Villard, Paris, 1893-1902.\\ 

[11] \ E.T. Whittaker,G.N. Watson \quad {\it A course of Modern Analysis}\\ Cambridge (1963).\\

\end{document}